\documentclass[reqno]{amsart}

\usepackage{graphicx}
\usepackage{amsmath}
\usepackage{amssymb}
\usepackage{amsthm}
\usepackage{enumitem}
\usepackage{bbm}
\usepackage[all,error]{onlyamsmath}
\usepackage[strict]{csquotes}
\usepackage{url}
\usepackage{braket,mleftright}
\usepackage[mathscr]{euscript}

\RequirePackage[final=true]{hyperref}
\hypersetup{colorlinks=true,allcolors=[rgb]{0,0,0.5}}

\newcommand{\ie}{\textit{i.e.}\;}
\newcommand{\eg}{\textit{e.g.}\;}
\newcommand{\cf}{\textit{cf.}\;}

\newcommand{\mrm}[1]{\mathrm{#1}}		
\newcommand{\ol}[1]{\overline{#1}}		
\newcommand{\vrt}{\,\vert\,}			
\newcommand{\what}[1]{\widehat{#1}}		
\newcommand{\wtilde}[1]{\widetilde{#1}}	
\newcommand{\mcap}{\cap}				
\newcommand{\op}{\oplus}				
\newcommand{\dsum}{\dotplus}			
\newcommand{\hsum}{\,\what{+}\,}		
\newcommand{\pd}{\partial}				
\DeclareMathOperator{\dom}{dom}		
\DeclareMathOperator{\ran}{ran}		
\DeclareMathOperator{\mul}{mul}		

\newcommand{\fH}{\mathfrak{H}}

\newcommand{\fK}{\mathfrak{K}}
\newcommand{\fL}{\mathfrak{L}}
\newcommand{\fM}{\mathfrak{M}}
\newcommand{\fN}{\mathfrak{N}}

\newcommand{\wtP}{\wtilde{P}}
\newcommand{\wtQ}{\wtilde{Q}}

\newcommand{\whh}{\what{h}}
\newcommand{\whk}{\what{k}}

\newcommand{\row}[2]{\begin{pmatrix}#1 & #2 \end{pmatrix} }
\newcommand{\col}[2]{\begin{pmatrix}#1 \\ #2 \end{pmatrix} }

\newtheorem{thm}{Theorem}
\newtheorem{lem}[thm]{Lemma}
\newtheorem{cor}[thm]{Corollary}

\theoremstyle{definition}
\newtheorem{defn}[thm]{Definition}

\newtheorem{exam}[thm]{Example}
\numberwithin{equation}{section}
\numberwithin{thm}{section}
\begin{document}

\title[On rows and columns of linear relations]{On necessary and
sufficient
conditions relating the adjoint of a column to a row of
linear relations}
\author{Rytis Jur\v{s}\.{e}nas}
\address{Vilnius University,
Institute of Theoretical Physics and Astronomy,
Saul\.{e}tekio ave.~3, LT-10257 Vilnius, Lithuania}
\email{Rytis.Jursenas@tfai.vu.lt}
\thanks{}
\date{\today}
\begin{abstract}
A row and a column of two linear relations in Hilbert spaces
are presented respectively as a sum and an intersection
of two linear relations. As an application, necessary and
sufficient conditions for the adjoint of
a column to be a row are examined. Several outcomes
are discussed as well.
\end{abstract}
\maketitle
\section{Introduction}
If $\col{C_1}{C_2}$ is a densely defined column of two
operators $C_1$ and $C_2$ then by
\cite[Proposition~4.3-3$^\circ$]{Moller08} the
adjoint of the column is the closure of the row
$\row{C^*_1}{C^*_2}$, \ie
\begin{equation}
\col{C_1}{C_2}^*=\ol{\row{C^*_1}{C^*_2}}
\label{eq:I1}
\end{equation}
iff the closure of the column
is the column of the closures $\ol{C_1}$ and $\ol{C_2}$;
we label the latter condition by \eqref{eq:C}.
Here the closure is in general assumed in the sense
of linear relations, so that $C_i$, $i\in\{1,2\}$, is allowed
to be nonclosable as an operator. In
\cite[Proposition~4.1]{Hassi20}, the statement
is extended to linear relations $C_i$, by showing that
the adjoint of the column is the row
$\row{C^*_1}{C^*_2}$, \ie
\begin{equation}
\col{C_1}{C_2}^*=\row{C^*_1}{C^*_2}
\label{eq:I2}
\end{equation}
iff the adjoint of the column is
\textit{a} row.

In this note we consider
columns and rows of two linear relations, but
it will become clear from the exposition that our method
extends naturally to the case of a finite number
$\geq2$ of linear relations.
We prove that
\eqref{eq:I2} is true iff \eqref{eq:C} holds
and in addition
$\dom\ol{C_1}+\dom\ol{C_2}$ is closed;
we label the latter condition by \eqref{eq:Cp}.
In other words the row $\row{C^*_1}{C^*_2}$ is closed
iff \eqref{eq:Cp} holds. Some consequences of the
theorem are also discussed. For example, it is known
from \cite[Corollary~3.4]{Moller08} that, for
a row $\row{R_1}{R_2}$ of two operators $R_1$ and $R_2$
to be closable it is necessary that $R_1$ and $R_2$
be closable. We obtain that, if in addition
$\dom R^*_1+\dom R^*_2$ is closed, then the condition
is also sufficient.

The main technical ingredient in the present
note is a demonstration that the column $\col{C_1}{C_2}$
appears to be the intersection $A_1\mcap A_2$
of linear relations $A_i:=Q^{-1}_iC_i$ for some
projections $Q_i$. Likewise, the row $\row{C^*_1}{C^*_2}$
is the componentwise sum $A^*_1\hsum A^*_2$, and
similarly for an arbitrary row $\row{R_1}{R_2}$.
From here we conclude that
\eqref{eq:I1} holds iff
$\ol{A_1\mcap A_2}=\ol{A_1}\mcap\ol{A_2}$, which
is \eqref{eq:C}. Moreover, because $A^*_1\hsum A^*_2$
is closed iff so is $\ol{A_1}\hsum\ol{A_2}$, we deduce
\eqref{eq:Cp}. The reader may refer \eg to
\cite{Jung19}
for the discussion as to when the equality
$\ol{A_1\mcap A_2}=\ol{A_1}\mcap\ol{A_2}$ holds true
for subsets of a topological space. For example,
it is true iff $\pd A_1\mcap\pd A_2\subseteq\pd(A_1\mcap A_2)$,
where $\pd$ denotes the boundary of a set.
However, we are not able to extract from this
a readable iff argument in full generality,
unless $A_1\subseteq A_2$ or $A_2\subseteq A_1$ or
$C_i$ (and hence $A_i$) is a singular linear relation.
On the other hand, if $\dom C_1\subseteq\dom C_2$,
$\ol{\mul}C_2=\mul\ol{C_2}$, and $\dom\ol{C_2}$ is closed,
then by \cite[Lemma~4.1]{Hassi20} these conditions
are sufficient to satisfy
\eqref{eq:C} and \eqref{eq:Cp}.

Let us mention that
rows and columns of operators or linear relations appear as building
blocks in matrix theory of operators, \eg
\cite{Moller08,Ota03,Nagel89}.
They prove to be useful in applications as well,
for example in
extension theory for sums of nonnegative linear relations
\cite[Lemma~3.1, Eq.~(3.8)]{Hassi07a}.
\section{Notation. Main tools}
Symbols $\fH$, $\fH_i$, $\fK$, $\fK_i$ denote Hilbert
spaces. The index $i\in\{1,2\}$. The scalar product in $\fH$
is denoted by $\braket{\cdot,\cdot}_\fH$.
The Cartesian product space $\fH\times\fK$ is identified
with the Hilbert sum $\fH\op\fK$ equipped with the usual
cross-product topology. A linear relation from
$\fH$ to $\fK$ is a linear subset of $\fH\times\fK$
with domain, range, kernel, multivalued part indicated
by $\dom$, $\ran$, $\ker$, $\mul$. The adjoint
linear relation is labeled by the asterisk, the
closure (double adjoint) by the overbar. A (linear) operator
is a linear relation with a trivial multivalued part.
An operator is closable if its closure is an operator.
The inverse $A^{-1}$ of a linear relation $A$ from
$\fH$ to $\fK$ is a linear relation from $\fK$ to $\fH$
which consists of pairs $(k,h)$ such that $(h,k)\in A$.
The inverse in the sense of linear relations exists
even if $\ker A$ is nontrivial. We always regard the inverses
in this way.

The (operatorwise) sum $A_1+A_2$ of linear relations
$A_i$ from $\fH$
to $\fK$ is a linear relation consisting of pairs
$(h,k_1+k_2)$ with $(h,k_i)\in A_i$. The
componentwise sum $A_1\hsum A_2$ consists of
$(h_1+h_2,k_1+k_2)$ with $(h_i,k_i)\in A_i$.
When computing the adjoint of a product and a
componentwise sum of
linear relations we rely on the next two lemmas:
\begin{lem}{\cite[Lemma~2.9]{Derkach09}}\label{lem:Derk}
Let $\fK_l$, $l\in\{0,1,2,3\}$, be Hilbert spaces,
$S\subseteq\fK_1\times\fK_2$ a closed
linear relation.
\begin{itemize}
\item[$\mrm{(i)}$]
If $\dom S$ is closed, then
$(SX)^*=X^*S^*$ for every $X\subseteq\fK_0\times\fK_1$
such that $\ran X\subseteq\dom S$.
\item[$\mrm{(ii)}$]
If $\ran S$ is closed, then
$(YS)^*=S^*Y^*$ for every $Y\subseteq\fK_2\times\fK_3$
such that $\dom Y\subseteq\ran S$.
\end{itemize}
\end{lem}
\begin{lem}{\cite[Lemma~2.10]{Hassi09}}\label{lem:Hassi}
Let $A_1$ and $A_2$ be linear relations in a Hilbert
space $\fH$.
The following statements are equivalent:
\begin{itemize}
\item[$\mrm{(i)}$]
$\ol{A_1}\hsum\ol{A_2}$ is closed;
\item[$\mrm{(ii)}$]
$A^*_1\hsum A^*_2$ is closed.
\end{itemize}
\end{lem}
We remark that Lemma~\ref{lem:Hassi} remains valid
for linear relations from $\fH$ to $\fK$.
Also
\begin{equation}
(A_1\hsum A_2)^*=A^*_1\mcap A^*_2
\quad\text{and hence}\quad
(\ol{A_1}\mcap\ol{A_2})^*=\ol{A^*_1\hsum A^*_2}
\label{eq:A1A2}
\end{equation}
for linear relations $A_i$ from $\fH$ to $\fK$.
\section{Main result}
We recall the definitions of the
main objects of our study.
\begin{defn}
Given $R_i\subseteq\fH_i\times\fK$,
the row of $R_1$ and $R_2$ is a linear relation
\[
\row{R_1}{R_2}\subseteq\fH\times\fK\,,\quad
\fH=\fH_1\times\fH_2
\]
defined by
\[
\row{R_1}{R_2}:=
\{((h_1,h_2),k_1+k_2)\in\fH\times\fK\vrt
(h_i,k_i)\in R_i\}\,.
\]
\end{defn}
\begin{defn}
Given $C_i\subseteq\fH\times\fK_i$,
the column of $C_1$ and $C_2$ is a linear relation
\[
\col{C_1}{C_2}\subseteq\fH\times\fK\,,\quad
\fK=\fK_1\times\fK_2
\]
defined by
\[
\col{C_1}{C_2}:=
\{(h,(k_1,k_2))\in\fH\times\fK\vrt
(h,k_i)\in C_i\}\,.
\]
\end{defn}
Next we consider the operator $P_i$ from
$\fH=\fH_1\times\fH_2$ to $\fH_i$
defined by $(h_1,h_2)\mapsto h_i$.
Notice that,
if we identify $\fH_1\times\{0\}$ with $\fH_1$
and $\{0\}\times\fH_2$ with $\fH_2$, then $P_i$
becomes an orthogonal projection in $\fH$ onto $\fH_i$.
However, in what follows
we regard $P_i$ as a single-valued linear relation from
$\fH$ to $\fH_i$. It is easy to verify that
in this case the adjoint operator
$P^*_i$ from $\fH_i$ to $\fH$ is described as follows:
\begin{lem}\label{lem:Pi}
The adjoint $P^*_i$ of $P_i$
is the operator from $\fH_i$ to $\fH$ given by
\[
P^*_1=I_{\fH_1}\op0\,,\quad
P^*_2=0\op I_{\fH_2}\,.
\]
In particular, $\wtP_i:=P^*_i$ is isometric, \ie
$\wtP^{-1}_i\subseteq\wtP^*_i$.
\end{lem}
Equivalently one may say that $P_i$ is coisometric.
\begin{proof}
As a linear relation, $P^*_1$
consists of $(h^\prime_1,\whh_\star)\in\fH_1\times\fH$,
with $\whh_\star=(h_\star,h^\prime_\star)$,
such that $(\forall\whh=(h_1,h_2)\in\fH)$
$\braket{\whh,\whh_\star}_\fH=
\braket{h_1,h^\prime_1}_{\fH_1}$; hence
$h_\star=h^\prime_1$ and $h^\prime_\star=0$.
Because $P_1P^*_1=I_{\fH_1}$ but
$P^*_1P_1=I_{\fH_1}\op0\subsetneq I_\fH$ we get that $P^*_1$ is
isometric.
Similar considerations apply to $P^*_2$.
\end{proof}
In case $\fH_i$ (resp. $\fH=\fH_1\times\fH_2$)
is replaced by $\fK_i$ (resp. $\fK=\fK_1\times\fK_2$)
we use the symbol $Q_i$ in place of $P_i$, as well as
$\wtQ_i$ in place of $\wtP_i$.
We do so because we keep the spaces
$\fK_i$ and $\fH_i$ fixed,
while $Q_i$ and $P_i$ act in (generally) different spaces.

We now describe a row and a column in terms
of standard operations of linear relations;
namely, a sum and an intersection.
\begin{lem}\label{lem:row2}
\[
\row{R_1}{R_2}=R_1\wtP^{-1}_1\hsum
R_2\wtP^{-1}_2=R_1P_1+R_2P_2\,.
\]
\end{lem}
\begin{proof}
We have
\[
\wtP^{-1}_1=\{((h_1,0),h_1)\vrt h_1\in\fH_1 \}\,,\quad
\wtP^{-1}_2=\{((0,h_2),h_2)\vrt h_2\in\fH_2 \}
\]
so
\begin{align*}
R_1\wtP^{-1}_1=&
\{((h_1,0),k_1)\vrt(h_1,k_1)\in R_1 \}\,,
\\
R_2\wtP^{-1}_2=&
\{((0,h_2),k_2)\vrt(h_2,k_2)\in R_2 \}
\end{align*}
and hence
\[
R_1\wtP^{-1}_1\hsum R_2\wtP^{-1}_2=
\row{R_1}{R_2}\,.
\]

By definition, $\row{R_1}{R_2}$ consists of
$(\whh,k_1+k_2)$ such that $(P_i\whh,k_i)\in R_i$.
But
\[
R_iP_i=\{(\whh,k_i)\vrt (P_i\whh,k_i)\in R_i \}
\]
so
\[
(P_i\whh,k_i)\in R_i\quad\Longleftrightarrow\quad
(\whh,k_i)\in R_iP_i
\]
and then
\begin{align*}
\row{R_1}{R_2}=&\{(\whh,k_1+k_2)\vrt
(\whh,k_i)\in R_iP_i\}
=R_1P_1+R_2P_2\,.\qedhere
\end{align*}
\end{proof}
\begin{lem}\label{lem:col}
\[
\col{C_1}{C_2}=(Q^{-1}_1C_1)\mcap(Q^{-1}_2C_2)\,.
\]
\end{lem}
\begin{proof}
By definition, $\col{C_1}{C_2}$ consists of
$(h,\whk)$ such that $(h,Q_i\whk)\in C_i$.
But
\[
Q^{-1}_iC_i=\{(h,\whk)\vrt (h,Q_i\whk)\in C_i \}
\]
so
\[
(h,Q_i\whk)\in C_i\quad\Longleftrightarrow\quad
(h,\whk)\in Q^{-1}_iC_i
\]
and then
\begin{align*}
\col{C_1}{C_2}=&\{(h,\whk)\vrt
(h,\whk)\in Q^{-1}_iC_i\}
=(Q^{-1}_1C_1)\mcap(Q^{-1}_2C_2)\,.\qedhere
\end{align*}
\end{proof}
\begin{thm}\label{thm:RC}
The following statements hold:
\begin{itemize}
\item[$\mrm{(i)}$]
$\row{R_1}{R_2}^*=\col{R^*_1}{R^*_2}$.
\medskip
\item[$\mrm{(ii)}$]
$\col{C_1}{C_2}^*\supseteq\ol{\row{C^*_1}{C^*_2}}$,
with the equality iff
\begin{equation}
\ol{\col{C_1}{C_2}}=\col{\ol{C_1}}{\ol{C_2}}\,.
\label{eq:C}\tag{$C$}
\end{equation}
\medskip
\item[$\mrm{(iii)}$]
$\row{C^*_1}{C^*_2}$ is closed
iff
\begin{equation}
\dom\ol{C_1}+\dom\ol{C_2}\quad\text{is closed}\,.
\label{eq:Cp}\tag{$C^\prime$}
\end{equation}
\end{itemize}
\end{thm}
\begin{proof}
(i)
By applying Lemma~\ref{lem:row2} and
\eqref{eq:A1A2}
\[
\row{R_1}{R_2}^*=(R_1\wtP^{-1}_1\hsum R_2\wtP^{-1}_2)^*=
(R_1\wtP^{-1}_1)^*\mcap(R_2\wtP^{-1}_2)^*\,.
\]
Because $\wtP^{-1}_i$ is closed, and because
$\dom \wtP_i=\fH_i$ is closed
and contains $\dom R_i$, by applying
Lemmas~\ref{lem:Derk}(ii) and \ref{lem:Pi}
we get that
\[
(R_i\wtP^{-1}_i)^*=\wtP^{*\,-1}_iR^*_i=
P^{-1}_iR^*_i\,.
\]
Therefore, by Lemma~\ref{lem:col}
\[
\row{R_1}{R_2}^*=(P^{-1}_1R^*_1)\mcap
(P^{-1}_2R^*_2)=\col{R^*_1}{R^*_2}\,.
\]

(ii)
Let $A_i:=Q^{-1}_iC_i$.
Observe that the adjoint and the closure
\[
A^*_i=(Q^{-1}_iC_i)^*=
C^*_i\wtQ^{-1}_i\quad\text{and}\quad
\ol{A_i}=\ol{Q^{-1}_iC_i}=Q^{-1}_i\ol{C_i}
\]
by applying Lemma~\ref{lem:Derk}.
We have by Lemma~\ref{lem:col} and \eqref{eq:A1A2}
\[
\col{C_1}{C_2}=
A_1\mcap A_2\subseteq
\ol{A_1}\mcap\ol{A_2}
=(A^*_1\hsum A^*_2)^*\,.
\]
Taking the adjoints this implies that
\[
\col{C_1}{C_2}^*=
(A_1\mcap A_2)^*
\supseteq
\ol{A^*_1\hsum A^*_2}\,.
\]
By applying Lemmas~\ref{lem:Pi} and \ref{lem:row2}
\[
\col{C_1}{C_2}^*\supseteq
\ol{C^*_1\wtQ^{-1}_1\hsum C^*_2\wtQ^{-1}_2}=
\ol{\row{C^*_1}{C^*_2}}\supseteq
\row{C^*_1}{C^*_2}\,.
\]

Next, because
\[
\col{C_1}{C_2}^*=(A_1\mcap A_2)^*=
(\ol{A_1\mcap A_2})^*\,,
\]
\[
(\ol{A_1}\mcap\ol{A_2})^*
=\ol{A^*_1\hsum A^*_2}=
\ol{\row{C^*_1}{C^*_2}}
\]
we have the following: If \eqref{eq:C}, that is,
if $\ol{A_1\mcap A_2}=\ol{A_1}\mcap\ol{A_2}$, then
\[
\col{C_1}{C_2}^*=\ol{\row{C^*_1}{C^*_2}}\,.
\]
Conversely, if the latter equality holds, then
\[
(A_1\mcap A_2)^*=(\ol{A_1}\mcap\ol{A_2})^*
\quad\Longrightarrow\quad
\ol{A_1\mcap A_2}=
\ol{\ol{A_1}\mcap\ol{A_2}}=\ol{A_1}\mcap\ol{A_2}
\]
\ie \eqref{eq:C} holds.

(iii)
Finally, $\row{C^*_1}{C^*_2}$ is closed iff
$A^*_1\hsum A^*_2$ is closed. According to
Lemma~\ref{lem:Hassi} this can happen iff
$\ol{A_1}\hsum\ol{A_2}$ is closed.
But
\[
\ol{A_1}\hsum\ol{A_2}=
(\dom\ol{C_1}+\dom\ol{C_2})\times\fK\,,\quad
\fK=\fK_1\times\fK_2
\]
so $\row{C^*_1}{C^*_2}$ is closed iff
\eqref{eq:Cp}.
\end{proof}
We mention that (i) could be also shown by using the
second equality in Lemma~\ref{lem:row2} and then
by applying Lemma~\ref{lem:Pi}.

As already remarked, (i) is known from
\cite[Proposition~2.1(i)]{Hassi07a}
(see also \cite[Proposition~4.1]{Hassi20});
(ii) is stated in \cite[Proposition~4.3-3$^\circ$]{Moller08}
for a densely defined column of operators;
(iii) seems to be new (to the best of our knowledge).
\section{Some corollaries}
In this section we discuss some consequences
following from Theorem~\ref{thm:RC}.
\begin{cor}\label{cor:col}
The following statements are equivalent:
\begin{itemize}
\item[$\mrm{(i)}$]
$\col{C_1}{C_2}^*=\row{C^*_1}{C^*_2}$.
\medskip
\item[$\mrm{(ii)}$]
\eqref{eq:C} and \eqref{eq:Cp} hold.
\medskip
\item[$\mrm{(iii)}$]
$\col{C_1}{C_2}^*$ is a row.
\end{itemize}
\end{cor}
\begin{proof}
(i) $\Leftrightarrow$ (ii) is due to Theorem~\ref{thm:RC}.
(i) $\Rightarrow$ (iii) is obvious.
(iii) $\Rightarrow$ (i) is shown in
\cite[Proposition~4.1]{Hassi20}.
\end{proof}
Let us repeat that \eqref{eq:C} is equivalent to
\[
\ol{A_1\mcap A_2}=\ol{A_1}\mcap\ol{A_2}\quad
\text{with}\quad
A_i:=Q^{-1}_iC_i\,,\quad
\ol{A_i}=Q^{-1}_i\ol{C_i}\,.
\]
This equality does not hold for all $C_i$,
as in general one only has the inclusion
$\subseteq$. Therefore, \eqref{eq:C} is equivalent to
showing that $\ol{A_1\mcap A_2}\supseteq\ol{A_1}\mcap\ol{A_2}$.
In general it might be not easy to verify the latter
inclusion, but in some special cases the condition could
be stated more explicitly. The most obvious one
is \eg $A_1\subseteq A_2$, that is,
$\dom C_1\subseteq\dom C_2$ and
$\ran (C_2\vrt_{\dom C_1})=\fK_2$.
Another example is
$\dom C_1\subseteq\dom C_2$ but
$\ran (C_2\vrt_{\dom C_1})\subseteq\fK_2$
and $C_2$ is singular:
\begin{exam}(\cf \cite[Corollary~4.4]{Hassi20})
Assume that $C_2$ is a singular linear relation
(\cite[Eq.~(3.3)]{Hassi07}) of the form
$C_2=\fM\times\fN$ for some $\fM\subseteq\fH$
and $\fN\subseteq\fK_2$; $\fM$, $\fN$ need not be
closed. Then
\[
A_1\mcap A_2=C_1\vrt_{\dom C_1\mcap\fM}\times\fN\,,\quad
\ol{A_1}\mcap\ol{A_2}=
\ol{C_1}\vrt_{\dom\ol{C_1}\mcap\ol{\fM}}\times\ol{\fN}
\]
and conditions \eqref{eq:C} and \eqref{eq:Cp}
read
\[
(C)\;\;
\ol{C_1\vrt_{\dom C_1\mcap\fM}}=
\ol{C_1}\vrt_{\dom\ol{C_1}\mcap\ol{\fM}}
\quad\text{and}\quad
(C^\prime)\;\;
\dom\ol{C_1}+\ol{\fM}\;\text{is closed.}
\]
For instance, if
$\dom C_1\subseteq\fM$, then
\[
\dom\ol{C_1}\subseteq\ol{\dom}C_1=\ol{\dom}\;\ol{C_1}
\subseteq\ol{\fM}
\]
and therefore (i) in Corollary~\ref{cor:col} holds:
\begin{align*}
\col{C_1}{\fM\times\fN}^*=&
A^*_1\hsum(\{0\}\times\fN^\bot)
\times\fM^\bot
=\row{C^*_1}{\fN^\bot\times\fM^\bot}
\\
=&A^*_1\hsum(\{0\}\times\fN^\bot)
\times\{0\}=\row{C^*_1}{\fN^\bot\times\{0\}}
\\
=&\col{C_1}{\fH\times\fN}^*
\end{align*}
where $\fM^\bot$ (resp. $\fN^\bot$) is the
orthogonal complement in $\fH$ (resp. $\fK_2$)
of $\fM$ (resp. $\fN$). We remark that the third
equality is due to the implication
$\dom C_1\subseteq\fM$ $\Rightarrow$
$\mul C^*_1\supseteq\fM^\bot$, and that
$\fH\times\fN$ in the last equality can be replaced
by $\dom C_1\times\fN$ due to the definition of the row.
\end{exam}
In \cite[Lemma~4.1]{Hassi20}
sufficient conditions
to have (i) in Corollary~\ref{cor:col} are given. Namely:
\begin{enumerate}[label=(\emph{\alph*}),ref=(\emph{\alph*})]
\item\label{item:a}
$\dom C_1\subseteq\dom C_2$;
\medskip
\item\label{item:b}
$\ol{\mul}C_2=\mul\ol{C_2}$;
\medskip
\item\label{item:c}
$\dom\ol{C_2}$ is closed.
\end{enumerate}
We see that \ref{item:a} and \ref{item:c} together
imply \eqref{eq:Cp}.
We also see that in the previous example
\ref{item:b} need not hold.

Another example taken from \cite{Hassi20}
is the following:
\begin{exam}
Let $B_i\in[\fH]$ (bounded), $B_i>0$,
$\ran B_1\mcap\ran B_2=\{0\}$. Taking $C_i:=B^{-1}_i=\ol{C_i}$
we see that $\ran B_1\dsum\ran B_2$ (direct sum)
cannot be closed, since the spectrum
$\sigma(B_i)\subseteq[0,\infty]$ and $B^*_i=B_i$,
and therefore $\ran B_i$ is not closed. Thus
\eqref{eq:Cp} does not hold. But
\eqref{eq:C} does hold, since $A_i=(B_iQ_i)^{-1}$ is closed;
hence
\[
\col{B^{-1}_1}{B^{-1}_2}^*=
\ol{\row{B^{*\,-1}_1}{B^{*\,-1}_2}}
\supsetneq\row{B^{*\,-1}_1}{B^{*\,-1}_2}\,.
\]
Note that in this example $\fK_1=\fK_2=\fH$.
\end{exam}
It follows from Theorem~\ref{thm:RC}(i) that
the column of two closed linear relations is closed.
In particular, the column of two closable operators
is a closable operator.
For a row we have a different situation
(\cf \cite[Corollaries~4.1, 4.2]{Hassi20}):
\begin{cor}\label{cor:R}
The following statements hold:
\begin{itemize}
\item[$\mrm{(i)}$]
$\ol{\row{R_1}{R_2}}=\ol{\row{\ol{R_1}}{\ol{R_2}}}$.
\medskip
\item[$\mrm{(ii)}$]
$\row{\ol{R_1}}{\ol{R_2}}\subseteq
\ol{\row{R_1}{R_2}}$, with the equality
iff
\begin{equation}
\dom R^*_1+\dom R^*_2\quad\text{is closed}\,.
\label{eq:R*}\tag{$R$}
\end{equation}
\end{itemize}
\end{cor}
\begin{proof}
(i)
Because $R^*_i$ is closed, the column of
$R^*_1$ and $R^*_2$ is closed; hence \eqref{eq:C}
holds (with $R^*_i$ in place of $C_i$).
Thus the equality follows from
Theorem~\ref{thm:RC}(i) and (ii).

(ii)
This is due to (i) and Theorem~\ref{thm:RC}(iii)
(with $R^*_i$ in place of $C_i$).
\end{proof}
In particular, the equality in Corollary~\ref{cor:R}(ii)
holds if
$\dom R^*_1\subseteq\dom R^*_2$ and $\dom R^*_2$
is closed.

Because
\[
\mul\row{R_1}{R_2}=\mul R_1+\mul R_2
\]
it follows that the row is an operator iff
$R_i$ is an operator for each $i\in\{1,2\}$.
In this direction we deduce the next corollary
(\cf \cite[Corollary~3.4]{Moller08}).
\begin{cor}\label{cor:row2}
Let $R_1$ and $R_2$ be operators. For a row
$\row{R_1}{R_2}$ to be closable it is necessary and,
in case \eqref{eq:R*} holds, also
sufficient that $R_1$ and $R_2$ are closable.
\end{cor}
\begin{proof}
Necessity:
Assume that $\row{R_1}{R_2}$ is closable, that is,
$\ol{\row{R_1}{R_2}}$ is an operator. Then
by Corollary~\ref{cor:R}(ii)
\[
\mul\row{\ol{R_1}}{\ol{R_2}}=
\mul \ol{R_1}+\mul \ol{R_2}\subseteq
\mul\ol{\row{R_2}{R_2}}=\{0\}\,;
\]
hence $\mul\ol{R_i}=\{0\}$ and $R_i$ is closable.

Sufficiency: Assume $\mul\ol{R_i}=\{0\}$ and \eqref{eq:R*}.
Then again by Corollary~\ref{cor:R}(ii), the closure
$\ol{\row{R_1}{R_2}}=\row{\ol{R_1}}{\ol{R_2}}$ is an operator.
\end{proof}
We note that, if $R_1$ and $R_2$ are closable operators and
\eqref{eq:R*} holds, then necessarily
$\dom R^*_1+\dom R^*_2=\fK$. To see this, let
$\fL_i$ be dense linear subsets of $\fK$, and such that
$\fL_1+\fL_2=:\fM$ is closed. Then
$\fM^\bot=\fL^\bot_1\mcap\fL^\bot_2=\{0\}$ implies that $\fM=\fK$.
If in particular $\dom R^*_1\subseteq\dom R^*_2$,
then $R_2$ is bounded everywhere defined, and we arrive
at \cite[Proposition~3.5]{Moller08}.

For not necessarily closed operators $C_1$ and $C_2$,
we trivially
have that the column of $C_1$ and $C_2$ is
closable iff $\ol{A_1\mcap A_2}$ is an operator.
On the other hand,
because
\[
\mul\col{C_1}{C_2}=\mul C_1\times\mul C_2
\]
to make $\col{C_1}{C_2}$ closable Theorem~\ref{thm:RC}(i)-(ii)
shows that it is necessary but not sufficient that $\ol{C_i}$
be an operator
(see also \cite[Proposition~4.3-2$^\circ$]{Moller08}).
Under additional conditions, however, the present
necessary condition becomes also sufficient.
\begin{cor}\label{cor:col2}
Let $C_1$ and $C_2$ be operators such that
$C_2$ is closable and \ref{item:a}, \ref{item:c} hold.
For a column
$\col{C_1}{C_2}$ to be closable it is necessary and
sufficient that $C_1$ is closable.
\end{cor}
\begin{proof}
By hypothesis, condition \eqref{eq:C} holds. Thus
\[
\mul\ol{\col{C_1}{C_2}}=\mul\ol{C_1}\times\{0\}
\]
and the claim follows.
\end{proof}
\section{Application to block relations}
As an example of application of Corollaries~\ref{cor:col},
\ref{cor:row2}, \ref{cor:col2}
to matrices of linear relations we have the next
two corollaries. We recall from
\cite{Hassi20} that a linear relation $A$ from
$\fH=\fH_1\times\fH_2$ to $\fK=\fK_1\times\fK_2$
generated by the block is defined by
\[
A:=\begin{pmatrix}
A_{11} & A_{12} \\ A_{21} & A_{22}
\end{pmatrix}=
\col{C_{A1}}{C_{A2}}=\row{R_{A1}}{R_{A2}}
\]
with
\[
C_{Ai}:=\row{A_{i1}}{A_{i2}}\,,\quad
R_{Ai}:=\col{A_{1i}}{A_{2i}}
\]
and $A_{ij}\subseteq\fH_j\times\fK_i$;
$i$, $j\in\{1,2\}$. Below we label by
$R^\prime_{Ai}$ the row
$\row{A^*_{1i}}{A^*_{2i}}$. Observe that
$R^\prime_{Ai}\subseteq R^*_{Ai}$ by Theorem~\ref{thm:RC}.
\begin{cor}
$1^\circ$ The adjoint
\begin{equation}
\begin{pmatrix}
A_{11} & A_{12} \\ A_{21} & A_{22}
\end{pmatrix}^*\supseteq
\begin{pmatrix}
A^*_{11} & A^*_{21} \\ A^*_{12} & A^*_{22}
\end{pmatrix}
\label{eq:A*}
\end{equation}
with the equality
iff
\begin{itemize}
\item[$(C_i)$]
$\ol{\col{A_{1i}}{A_{2i}}}=\col{\ol{A_{1i}}}{\ol{A_{2i}}}$
and
\medskip
\item[$(C^\prime_i)$]
$\dom\ol{A_{1i}}+\dom\ol{A_{2i}}$ is closed.
\end{itemize}

$2^\circ$
The closure
\[
\ol{\begin{pmatrix}
A_{11} & A_{12} \\ A_{21} & A_{22}
\end{pmatrix}}=
\begin{pmatrix}
\ol{A_{11}} & \ol{A_{12}} \\ \ol{A_{21}} & \ol{A_{22}}
\end{pmatrix}
\]
iff $(C_i)$ and $(C^\prime_i)$ hold and in addition
\begin{itemize}
\item[$(C^{\prime\prime}_i)$]
$\dom A^*_{i1}+\dom A^*_{i2}$ is closed.
\end{itemize}
\end{cor}
\begin{proof}
$1^\circ$ First note that the inclusion in
\eqref{eq:A*} is due to
\[
(A_{ij})^*=\col{C_{A1}}{C_{A2}}^*\supseteq
\row{C^*_{A1}}{C^*_{A2}}=(A^*_{ji})
\]
by Theorem~\ref{thm:RC}. According to
\cite[Corollary~6.1]{Hassi20}, there is equality
in \eqref{eq:A*} iff $R^*_{Ai}=R^\prime_{Ai}$,
which is $(C_i)$ and $(C^\prime_i)$ by
Corollary~\ref{cor:col}.

$2^\circ$ By $1^\circ$ we have
$\ol{(A_{ij})}=(\ol{A_{ij}})$ iff
$R^*_{Ai}=R^\prime_{Ai}$ and $R^*_{A^\prime i}=R^\prime_{A^\prime i}$,
with $A^\prime=(A^\prime_{ij})$
and $A^\prime_{ij}:=A^*_{ji}$; hence iff
$(C_i)$, $(C^\prime_i)$, and $(C^{\prime\prime}_i)$ hold.
\end{proof}
\begin{cor}
Let $A$ be an operator.
For $A$ to be closable
it is necessary and, in case
\[
\dom R^*_{A1}+\dom R^*_{A2}\quad\text{is closed},
\]
also sufficient that $R_{A1}$ and $R_{A2}$ are closable.
In particular, if
\begin{itemize}
\item[$(a_i)$]
$\dom A_{1i}\subseteq\dom A_{2i}$ and
\item[$(b_i)$]
$A_{2i}$ is closable and
\item[$(c_i)$]
$\dom\ol{A_{2i}}$ is closed
\end{itemize}
then $R_{Ai}$ is closable iff so is $A_{1i}$.
\end{cor}
\begin{proof}
This is an application of Corollaries~\ref{cor:row2},
\ref{cor:col2}
to $\row{R_{A1}}{R_{A2}}$.
\end{proof}


\end{document}